\documentclass{amsart}

\usepackage{amssymb,amsmath,amsthm,latexsym}

  \theoremstyle{definition}
  \newtheorem{definition}             {Definition}[section]
  
  \newtheorem{remark}     [definition]{Remark}

  \theoremstyle{plain}
  \newtheorem{lemma}      [definition]{Lemma}
  \newtheorem{proposition}[definition]{Proposition}
  \newtheorem{theorem}    [definition]{Theorem}
  \newtheorem{corollary}  [definition]{Corollary}

\begin{document}

\title[The nilpotent Malcev algebra of dimension five]
{Enveloping algebras of the nilpotent Malcev algebra of dimension five}

\author{Murray R. Bremner}

\address{Department of Mathematics and Statistics, University of Saskatchewan,
Saskatoon, Canada}

\email{bremner@math.usask.ca}

\author{Hamid Usefi}

\address{Department of Mathematics, University of British Columbia, Vancouver, Canada}

\email{usefi@math.ubc.ca}



\maketitle

\begin{abstract}
P\'erez-Izquierdo and Shestakov recently extended the PBW theorem to Malcev
algebras. It follows from their construction that for any Malcev algebra $M$
over a field of characteristic $\ne 2, 3$ there is a representation of the
universal nonassociative enveloping algebra $U(M)$ by linear operators on the
polynomial algebra $P(M)$. For the nilpotent non-Lie Malcev algebra
$\mathbb{M}$ of dimension 5, we use this representation to determine explicit
structure constants for $U(\mathbb{M})$; from this it follows that
$U(\mathbb{M})$ is not power-associative. We obtain a finite set of generators
for the alternator ideal $I(\mathbb{M}) \subset U(\mathbb{M})$ and derive
structure constants for the universal alternative enveloping algebra
$A(\mathbb{M}) = U(\mathbb{M})/I(\mathbb{M})$, a new infinite dimensional
alternative algebra. We verify that the map $\iota\colon \mathbb{M} \to
A(\mathbb{M})$ is injective, and so $\mathbb{M}$ is special.
\end{abstract}


\section{Introduction}

A nonassociative algebra $A$ is alternative if it satisfies $(x,x,y) = 0$ and
$(y,x,x) = 0$ for all $x, y \in A$ where $(x,y,z) = (xy)z - x(yz)$. Equivalent
conditions are that any subalgebra generated by two elements is associative,
and that the associator $(x,y,z)$ is an alternating function of its arguments.
If $A$ is an alternative algebra, then $A^-$ denotes the commutator algebra:
the same vector space with the operation $[x,y] = xy - yx$. The algebra $A^-$
satisfies anticommutativity $[x,x] = 0$ and the Malcev identity $[J(x,y,z),x] =
J(x,y,[x,z])$ where $J(x,y,z) = [[x,y],z] + [[y,z],x] + [[z,x],y]$; these two
identities define Malcev algebras. Basic references on Malcev algebras are
\cite{BremnerMurakamiShestakov}, \cite{Kuzmin1}, \cite{Malcev}, \cite{Sagle}.
The major unsolved problem in the theory of Malcev algebras is to determine
whether every Malcev algebra is isomorphic to a subalgebra of $A^-$ for some
alternative algebra $A$; that is, whether every Malcev algebra is special.
Substantial progress was made recently by P\'erez-Izquierdo and Shestakov
\cite{PerezIzquierdoShestakov}: they extended the Poincar\'e--Birkhoff--Witt
theorem to Malcev algebras by constructing universal nonassociative enveloping
algebras for Malcev algebras.

The smallest nilpotent non-Lie Malcev algebra $\mathbb{M}$ over a field of
characteristic 0 has dimension 5 and is unique up to isomorphism
\cite{Kuzmin2}. Its structure constants, which define a nilpotent Malcev
algebra over any field, are given in Table \ref{structureconstants}. In this
paper we construct a representation of the universal nonassociative enveloping
algebra $U(\mathbb{M})$ by differential operators on the polynomial algebra
$P(\mathbb{M})$; from this we determine explicit structure constants for
$U(\mathbb{M})$. We then determine a set of generators for the alternator ideal
$I(\mathbb{M})$; from this we obtain explicit structure constants for the
universal alternative enveloping algebra $A(\mathbb{M}) =
U(\mathbb{M})/I(\mathbb{M})$, and verify the speciality of $\mathbb{M}$.

  \begin{table}
  \begin{center}
  \label{structureconstants}
  \caption[]{The nilpotent Malcev algebra $\mathbb{M}$ of dimension five}
  \begin{tabular}{r|rrrrr}
  $[-{,}-]$ &\phantom{-} $\phantom{-}a$ & $\phantom{-}b$
  & $\phantom{-}c$ & $\phantom{-}d$ & $\phantom{-}e$
  \\
  \hline
  &&&&& \\[-8pt]
  $a$   & $ 0$ & $ c$ & $ 0$ & $ 0$ & $ 0$ \\
  $b$   & $-c$ & $ 0$ & $ 0$ & $ 0$ & $ 0$ \\
  $c$   & $ 0$ & $ 0$ & $ 0$ & $ e$ & $ 0$ \\
  $d$   & $ 0$ & $ 0$ & $-e$ & $ 0$ & $ 0$ \\
  $e$   & $ 0$ & $ 0$ & $ 0$ & $ 0$ & $ 0$
  \end{tabular}
  \end{center}
  \end{table}


\section{Theorem of P\'erez-Izquierdo and Shestakov}

All multilinear structures are over a field $\mathbb{F}$ of characteristic $\ne
2, 3$.

\begin{definition}
The \textbf{generalized alternative nucleus} of an algebra $A$ is
  \[
  N_{\mathrm{alt}}(A) =
  \big\{ \,
  a \in A \,
  | \,
  (a,x,y) = -(x,a,y) = (x,y,a), \,
  \text{for all} \,
  x, y \in A \,
  \big\}.
  \]
This is a subalgebra of $A^-$ (but not of $A$) and is a Malcev algebra.
\end{definition}

\begin{theorem}
[P\'erez-Izquierdo and Shestakov \cite{PerezIzquierdoShestakov}]
\label{PStheorem} For every Malcev algebra $M$ over $\mathbb{F}$ there exists a
nonassociative algebra $U(M)$ and an injective homomorphism $\iota\colon M \to
U(M)^-$ such that $\iota(M) \subseteq N_{\mathrm{alt}}(U(M))$; furthermore,
$U(M)$ is a universal object with respect to such homomorphisms.
\end{theorem}

Let $F(M)$ be the unital free nonassociative algebra over $\mathbb{F}$ on a
basis of $M$.  Let $R(M)$ be the ideal of $F(M)$ generated by the relations
  \[
  ab - ba - [a,b],
  \qquad
  (a,x,y) + (x,a,y),
  \qquad
  (x,a,y) + (x,y,a),
  \]
for all $a, b \in M$ and all $x, y \in F(M)$. Define $U(M) = F(M)/R(M)$ with
  \[
  \iota\colon M \to N_{\mathrm{alt}}(U(M)) \subseteq U(M),
  \qquad
  a \mapsto \iota(a) = \overline{a} = a + R(M).
  \]
Since $\iota$ is injective, we identify $M$ with $\iota(M) \subseteq U(M)$. Let
$B = \{ a_i \,|\, i \in \mathcal{I} \}$ be a basis of $M$ and let $<$ be a
total order on $\mathcal{I}$. Define
  \[
  \Omega =
  \{ \,
  (i_1,\hdots,i_n) \,
  | \,
  n \ge 0; \,
  i_1 \le \cdots \le i_n; \,
  i_1, \hdots, i_n \in \mathcal{I} \,
  \}.
  \]
For $n = 0$ the empty $n$-tuple gives $\overline{a}_\emptyset = 1 \in U(M)$.
For $n \ge 1$ the $n$-tuple $I = (i_1,\hdots,i_n) \in \Omega$ defines a
left-tapped monomial
  \[
  \overline{a}_I =
  \overline{a}_{i_1} (
  \overline{a}_{i_2} ( \cdots
  ( \overline{a}_{i_{n-1}} \overline{a}_{i_n} )
  \cdots )),
  \qquad
  |\overline{a}_I| = n.
  \]
The set of all $\overline{a}_I$ for $I \in \Omega$ is a basis of $U(M)$. For
any $f, g \in M$, $y \in U(M)$ we write an associator using commutators:
  \[
  (f,g,y)
  =
  \tfrac16
  [[y,f],g]
  -
  \tfrac16
  [[y,g],f]
  -
  \tfrac16
  [[y,[f,g]].
  \]
The next three Lemmas, which are implicit in \cite{PerezIzquierdoShestakov},
follow from the last formula and show inductively how to multiply in $U(M)$. We
first compute $[x,f]$ in $U(M)$; for $|x| = 1$ we use the product in $M$.

\begin{lemma}\label{rightbracket}
Let $x$ be a basis monomial of $U(M)$ with $|x| \ge 2$; write $x = gy$ with $g
\in M$. For any $f \in M$ we have
  \[
  [x,f]
  =
  [g,f]y
  +
  g[y,f]
  +
  \tfrac12
  [[y,f],g]
  -
  \tfrac12
  [[y,g],f]
  -
  \tfrac12
  [y,[f,g]].
  \]
\end{lemma}

We next compute $fx$ in $U(M)$; for $|x| = 1$ we have two cases: if $f \le x$
in the ordered basis of $M$, then $fx$ is an ordered monomial; otherwise, $fx =
xf + [f,x]$ where $[f,x] \in M$.

\begin{lemma}\label{leftmultiplication}
Let $x$ be a basis monomial of $U(M)$ with $|x| \ge 2$; write $x = gy$ with $g
\in M$. For any $f \in M$ we have
  \[
  fx =
  g(fy)
  +
  [f,g]y
  -
  \tfrac13
  [[y,f],g]
  +
  \tfrac13
  [[y,g],f]
  +
  \tfrac13
  [y,[f,g]].
  \]
\end{lemma}

We finally compute $yz$ in $U(M)$; for $|y| = 1$ we use Lemma
\ref{leftmultiplication}.

\begin{lemma}\label{generalproduct}
Let $y$ and $z$ be basis monomials of $U(M)$ with $|y| \ge 2$; write $y = fx$
with $f \in M$. We have
  \[
  yz = 2 f(xz) - x(fz) - x[z,f] + [xz,f].
  \]
\end{lemma}


\section{Representation by differential operators}

We write $P(M)$ for the polynomial algebra on the vector space $M$. Theorem
\ref{PStheorem} gives a linear isomorphism
  \[
  \phi\colon U(M) \to P(M),
  \quad
  \overline{a}_{i_1}
  (
  \overline{a}_{i_2}
  ( \cdots
  (
  \overline{a}_{i_{n-1}} \overline{a}_{i_n}
  )
  \cdots ))
  \mapsto
  a_{i_1}
  a_{i_2}
  \cdots
  a_{i_{n-1}} a_{i_n}.
  \]
In what follows we identify $U(M)$ with $P(M)$ by means of the linear
isomorphism $\phi$. This allows us to write monomials in $U(M)$ without
parentheses: $x$ represents $\phi^{-1}(x)$.

\begin{definition}
We define \textbf{right bracket} and \textbf{left multiplication}
maps
  \[
  \rho\colon U(M) \to \mathrm{End}_\mathbb{F} P(M),
  \qquad
  L\colon U(M) \to \mathrm{End}_\mathbb{F} P(M).
  \]
For $x \in U(M)$ we write $\rho(x)$ and $L(x)$ for the linear operators on
$P(M)$ induced by $y \mapsto [y,x]$ and $y \mapsto xy$ in $U(M)$:
  \[
  \rho(x)(f) = \phi \big( [ \phi^{-1}(f), x ] \big)
  \; \text{for} \;
  f \in P(M);
  \quad
  L(x)(f) = \phi \big( x \phi^{-1}(f) \big)
  \; \text{for} \;
  f \in P(M).
  \]
\end{definition}

For the rest of this paper $\mathbb{M}$ is the Malcev algebra of Table
\ref{structureconstants} over $\mathbb{F}$. We will show how to represent
$\rho(x)$ and $L(x)$ as differential operators on $P(\mathbb{M})$. Throughout
we assume the linear order $a < b < c < d < e$ on the basis of $\mathbb{M}$.

\begin{definition}
For $x \in \{a,b,c,d,e\}$ we write $M_x$ for multiplication by $x$ in
$P(\mathbb{M})$ and $D_x$ for differentiation with respect to $x$ in
$P(\mathbb{M})$. The next lemma is immediate.
\end{definition}

\begin{lemma} \label{DMcommutators}
$[ D_x, D_y ] = 0, \; [ M_x, M_y ] = 0, \; [ D_x, M_y ] = 0 \, (x \ne y), \; [
D_x, M_x ] = 1$.
\end{lemma}

The set $\{c,d,e\}$ spans a nilpotent Lie subalgebra $\mathbb{N} \subset
\mathbb{M}$. It follows from \cite{PerezIzquierdoShestakov} that $\mathbb{N}$
generates a subalgebra of $U(\mathbb{M})$ isomorphic to its associative
universal enveloping algebra $U(\mathbb{N})$.

\begin{proposition} \label{nilpotentLiealgebra}
In the associative subalgebra $U(\mathbb{N})$ of $U(\mathbb{M})$ we have
  \[
  ( c^i d^j e^k ) ( c^\ell d^m e^n )
  =
  \sum_{\alpha=0}^{\min(j,\ell)}
  (-1)^\alpha
  \alpha!
  \binom{j}{\alpha}
  \binom{\ell}{\alpha}
  c^{i+\ell-\alpha} d^{j+m-\alpha} e^{k+n+\alpha}.
  \]
\end{proposition}

\begin{proof}
We use a differential operator to illustrate the methods we apply later to the
nonassociative case. We first show that $d^j c = \big( M_c - M_e D_d \big) d^j$
by induction on $j$; the basis $j = 0$ is trivial. We use $[c,d] = e$, $[d,e] =
0$ to get
  \[
  d^{j+1} c
  =
  d ( c d^j {-} j d^{j-1} e )
  =
  c d^{j+1} - (j{+}1) d^j e
  =
  \big( M_c - M_e D_d \big) d^{j+1}.
  \]
We next show that $d^j c^\ell = \big( M_c {-} M_e D_d \big)^\ell d^j$ by
induction on $\ell$; the basis $\ell = 0$ is trivial. We use $[c,e] = 0$ to get
  \[
  d^j c^{\ell+1}
  =
  \big( ( M_c {-} M_e D_d )^\ell d^j \big) c
  =
  \big( M_c {-} M_e D_d \big)^\ell ( d^j c )
  =
  \big( M_c {-} M_e D_d \big)^{\ell+1} d^j.
  \]
We use $[c,e] = [d,e] = 0$ to get $( c^i d^j e^k ) ( c^\ell d^m e^n ) = c^i (
d^j c^\ell ) d^m e^{k+n}$. We now apply the binomial theorem since the terms
$M_c$ and $M_e D_d$ commute.  \end{proof}

\begin{corollary} \label{nilpotentcorollary}
We have $d ( b^q c^r d^s e^t ) = b^q c^r d^{s+1} e^t - r b^q c^{r-1}
d^s e^{t+1}$, and
  \begin{alignat*}{2}
  [ b^q c^r d^s e^t, b ] &= 0, &\quad
  [ b^q c^r d^s e^t, c ] &= - s b^q c^r d^{s-1} e^{t+1},
  \\
  [ b^q c^r d^s e^t, d ] &= r b^q c^{r-1} d^s e^{t+1}, &\quad
  [ b^q c^r d^s e^t, e ] &= 0,
  \\
  a ( b^q c^r d^s e^t ) &= a b^q c^r d^s e^t, &\quad
  b ( b^q c^r d^s e^t ) &= b^{q+1} c^r d^s e^t,
  \\
  c ( b^q c^r d^s e^t ) &= b^q c^{r+1} d^s e^t, &\quad
  e ( b^q c^r d^s e^t ) &= b^q c^r d^s e^{t+1}.
  \end{alignat*}
\end{corollary}

\begin{proof}
Since $b,c,d,e$ span a Lie subalgebra of $\mathbb{M}$, and $b$ commutes with
$c,d,e$, these all follow from Proposition \ref{nilpotentLiealgebra}.
\end{proof}

\begin{lemma} \label{bcde-brackets}
We have
  $\displaystyle{
  [ b^q c^r d^s e^t, a ]
  =
  -
  q
  b^{q-1} c^{r+1} d^s e^t
  +
  \tfrac{qs}{2}
  b^{q-1} c^r d^{s-1} e^{t+1}
  }$.
\end{lemma}

\begin{proof}
Lemma \ref{rightbracket} with $f = a$ and Corollary
\ref{nilpotentcorollary} give
  \allowdisplaybreaks
  \begin{align*}
  &
  [ b^{q+1} c^r d^s e^t, a ]
  =
  -
  c ( b^q c^r d^s e^t )
  +
  b [ b^q c^r d^s e^t, a ]
  \\
  &\qquad\qquad\qquad\qquad
  +
  \tfrac12
  [ [ b^q c^r d^s e^t, a ], b ]
  -
  \tfrac12
  [ [ b^q c^r d^s e^t, b ], a ]
  -
  \tfrac12
  [ b^q c^r d^s e^t, c ]
  \\
  &=
  -
  b^q c^{r+1} d^s e^t
  +
  b [ b^q c^r d^s e^t, a ]
  +
  \tfrac12
  [ [ b^q c^r d^s e^t, a ], b ]
  +
  \tfrac{s}{2}
  b^q c^r d^{s-1} e^{t+1}
  \\
  &=
  -
  b^q c^{r+1} d^s e^t
  +
  b
  \left(
  {-} q b^{q-1} c^{r+1} d^s e^t {+} \tfrac{qs}{2} b^{q-1} c^r d^{s-1} e^{t+1}
  \right)
  +
  \tfrac{s}{2}
  b^q c^r d^{s-1} e^{t+1}
  \\
  &
  =
  -
  (q{+}1)
  b^q c^{r+1} d^s e^t
  +
  \tfrac{(q{+}1)s}{2}
  b^q c^r d^{s-1} e^{t+1},
  \end{align*}
which completes the induction.  \end{proof}

\begin{proposition} \label{abcde-brackets}
As operators on $P(\mathbb{M})$ we have
  \allowdisplaybreaks
  \begin{alignat*}{2}
  (1) \;
  \rho(a) &= - M_c D_b + \tfrac12 M_e D_b D_d,
  &\quad
  (2) \;
  L(a) &= M_a,
  \\
  (3) \;
  \rho(b) &= M_c D_a - \tfrac12 M_e D_a D_d,
  &\quad
  (4) \;
  L(b) &= M_b - M_c D_a + \tfrac13 M_e D_a D_d,
  \\
  (5) \;
  \rho(c) &= - M_e D_d,
  &\quad
  (6) \;
  L(c) &= M_c,
  \\
  (7) \;
  \rho(d) &= M_e D_c + \tfrac12 M_e D_a D_b,
  &\quad
  (8) \;
  L(d) &= M_d - M_e D_c - \tfrac13 M_e D_a D_b,
  \\
  (9) \;
  \rho(e) &= 0,
  &\quad
  (10) \;
  L(e) &= M_e.
  \end{alignat*}
\end{proposition}

\begin{proof}
We use induction on $p$, the exponent of $a$; the basis $p = 0$ is Corollary
\ref{nilpotentcorollary} and Lemma \ref{bcde-brackets}. For (1), Lemma
\ref{rightbracket} with $f = a$ gives
  \[
  [ a^{p+1} b^q c^r d^s e^t, a ]
  =
  [ a ( a^p b^q c^r d^s e^t ), a ]
  =
  a [ a^p b^q c^r d^s e^t, a ].
  \]
Equation (2) is trivial. For (5), Lemma \ref{rightbracket} with $f =
c$ gives
  \allowdisplaybreaks
  \begin{align*}
  [ a^{p+1} b^q c^r d^s e^t, c ]
  &=
  a [ a^p b^q c^r d^s e^t, c ]
  {+}
  \tfrac12
  [ [ a^p b^q c^r d^s e^t, c ], a ]
  {-}
  \tfrac12
  [ [ a^p b^q c^r d^s e^t, a ], c ].
  \end{align*}
Applying (1) and induction we get
  \allowdisplaybreaks
  \begin{align*}
  &
  \rho(c)(a^{p+1} b^q c^r d^s e^t)
  =
  - M_a ( M_e D_d ) ( a^p b^q c^r d^s e^t )
  \\
  &\quad
  +
  \tfrac12
  ( M_c D_b )
  ( M_e D_d )
  ( a^p b^q c^r d^s e^t )
  -
  \tfrac14
  ( M_e D_b D_d )
  ( M_e D_d )
  ( a^p b^q c^r d^s e^t )
  \\
  &\quad\quad
  -
  \tfrac12
  ( M_e D_d )
  ( M_c D_b )
  ( a^p b^q c^r d^s e^t )
  +
  \tfrac14
  ( M_e D_d )
  ( M_e D_b D_d )
  ( a^p b^q c^r d^s e^t ).
  \end{align*}
By Lemma \ref{DMcommutators} the last four terms cancel and we get
  \[
  \rho(c)(a^{p+1} b^q c^r d^s e^t)
  =
  -
  M_a M_e D_d ( a^p b^q c^r d^s e^t )
  =
  -
  M_e D_d ( a^{p+1} b^q c^r d^s e^t ).
  \]
For (6), Lemma \ref{leftmultiplication} with $f = c$ gives
  \allowdisplaybreaks
  \begin{align*}
  c ( a^{p+1} b^q c^r d^s e^t )
  &=
  a ( c ( a^p b^q c^r d^s e^t ) )
  {-}
  \tfrac13
  [ [ a^p b^q c^r d^s e^t, c ], a ]
  {+}
  \tfrac13
  [ [ a^p b^q c^r d^s e^t, a ], c ].
  \end{align*}
Applying (1), (2), (5) and induction gives
  \allowdisplaybreaks
  \begin{align*}
  &
  L(c)( a^{p+1} b^q c^r d^s e^t )
  =
  ( M_a M_c )
  ( a^p b^q c^r d^s e^t )
  \\
  &\quad
  -
  \tfrac13
  ( M_c D_b )
  ( M_e D_d )
  ( a^p b^q c^r d^s e^t )
  +
  \tfrac16
  ( M_e D_b D_d )
  ( M_e D_d )
  ( a^p b^q c^r d^s e^t )
  \\
  &\quad\quad
  +
  \tfrac13
  ( M_e D_d )
  ( M_c D_b )
  ( a^p b^q c^r d^s e^t )
  -
  \tfrac16
  ( M_e D_d )
  ( M_e D_b D_d )
  ( a^p b^q c^r d^s e^t ).
  \end{align*}
By Lemma \ref{DMcommutators} the last four terms cancel and we get
  \[
  L(c)( a^{p+1} b^q c^r d^s e^t )
  =
  M_a M_c ( a^p b^q c^r d^s e^t )
  =
  M_c ( a^{p+1} b^q c^r d^s e^t ).
  \]
For (3), Lemma \ref{rightbracket} with $f = b$ gives
  \allowdisplaybreaks
  \begin{align*}
  &
  [ a^{p+1} b^q c^r d^s e^t, b ]
  =
  c ( a^p b^q c^r d^s e^t )
  +
  a [ a^p b^q c^r d^s e^t, b ]
  \\
  &\qquad
  +
  \tfrac12
  [ [ a^p b^q c^r d^s e^t, b ], a ]
  -
  \tfrac12
  [ [ a^p b^q c^r d^s e^t, a ], b ]
  +
  \tfrac12
  [ a^p b^q c^r d^s e^t, c ].
  \end{align*}
Applying (1), (2), (5), (6) and induction gives
  \[
  \rho(b)(a^{p+1} b^q c^r d^s e^t)
  =
  \left(
  M_c {+} M_a M_c D_a {-} \tfrac12 M_a M_e D_a D_d {-} \tfrac12 M_e D_d
  \right)
  ( a^p b^q c^r d^s e^t ),
  \]
using Lemma \ref{DMcommutators}. We now observe that
  \allowdisplaybreaks
  \begin{align*}
  ( M_c + M_c M_a D_a )
  ( a^p b^q c^r d^s e^t )
  &=
  ( M_c D_a )
  ( a^{p+1} b^q c^r d^s e^t ),
  \\
  ( M_a M_e D_a D_d + M_e D_d )
  ( a^p b^q c^r d^s e^t )
  &=
  ( M_e D_a D_d )
  ( a^{p+1} b^q c^r d^s e^t ).
  \end{align*}
For (4), Lemma \ref{leftmultiplication} with $f = b$ gives
  \allowdisplaybreaks
  \begin{align*}
  &
  b ( a^{p+1} b^q c^r d^s e^t )
  =
  a ( b ( a^p b^q c^r d^s e^t ) )
  -
  c ( a^p b^q c^r d^s e^t )
  \\
  &\qquad
  -
  \tfrac13
  [ [ a^p b^q c^r d^s e^t, b ], a ]
  +
  \tfrac13
  [ [ a^p b^q c^r d^s e^t, a ], b ]
  -
  \tfrac13
  [ a^p b^q c^r d^s e^t, c ].
  \end{align*}
Applying (1), (2), (3), (5), (6) and induction gives
  \begin{align*}
  &
  L(b) ( a^{p+1} b^q c^r d^s e^t )
  =
  \\
  &\quad
  \left(
  M_a M_b - M_a M_c D_a + \tfrac13 M_a M_e D_a D_d - M_c + \tfrac13 M_e D_d
  \right)
  ( a^p b^q c^r d^s e^t ),
  \end{align*}
using Lemma \ref{DMcommutators}, and now
  \[
  M_a M_b = M_b M_a,
  \;
  M_a M_c D_a + M_c = M_c D_a M_a,
  \;
  M_a M_e D_a D_d + M_e D_d = M_e D_a D_d M_a.
  \]
The proofs of (7)--(10) are similar.  \end{proof}

\begin{corollary} \label{operatorcommutators}
The nonzero commutators of $L(x)$ and $\rho(x)$ are
  \allowdisplaybreaks
  \begin{alignat*}{2}
  [ L(a), L(b) ] &= M_c - \tfrac13 M_e D_d,
  &\qquad
  [ L(a), L(d) ] &= \tfrac13 M_e D_b,
  \\
  [ L(b), L(d) ] &= -\tfrac13 M_e D_a,
  &\qquad
  [ L(c), L(d) ] &= M_e,
  \\
  [ \rho(a), \rho(d) ] &= M_e D_b,
  &\qquad
  [ \rho(b), \rho(d) ] &= - M_e D_a,
  \\
  [ L(a), \rho(b) ] &= - M_c + \tfrac12 M_e D_d,
  &\qquad
  [ L(a), \rho(d) ] &= - \tfrac12 M_e D_b,
  \\
  [ L(b), \rho(a) ] &= M_c - \tfrac12 M_e D_d,
  &\qquad
  [ L(b), \rho(d) ] &= \tfrac12 M_e D_a,
  \\
  [ L(c), \rho(d) ] &= - M_e,
  &\qquad
  [ L(d), \rho(a) ] &= \tfrac12 M_e D_b,
  \\
  [ L(d), \rho(b) ] &= - \tfrac12 M_e D_a,
  &\qquad
  [ L(d), \rho(c) ] &= M_e.
  \end{alignat*}
\end{corollary}


\section{Multiplication of basis monomials}

For operators $D, E$ with $[[D,E],E] = 0$ we have $[ D, E^n ] = n E^{n-1} [ D,
E ]$. The next result follows immediately from associativity for operators on
$P(\mathbb{N})$, but is nontrivial for operators on $P(\mathbb{M})$.

\begin{proposition} \label{cde-product}
We have $L( c^k d^\ell e^m ) = L(c)^k L(d)^\ell L(e)^m$.
\end{proposition}

\begin{proof}
We first prove $L(e^m) = L(e)^m$ by induction on $m$; the basis $m = 0$ is
trivial. We write $Z = a^p b^q c^r d^s e^t$. Lemma \ref{generalproduct} gives
  \[
  ( e^{m+1} ) Z
  =
  2 e ( e^m Z )
  -
  e^m ( e Z )
  -
  e^m [ Z, e ]
  +
  [ e^m Z, e ].
  \]
We now apply Proposition \ref{abcde-brackets}. We next prove $L( d^\ell e^m ) =
L(d)^\ell L(e)^m$ by induction on $\ell$; we have just proved the basis $\ell =
0$. Lemma \ref{generalproduct} gives
  \[
  ( d^{\ell+1} e^m ) Z
  =
  2
  d ( ( d^\ell e^m ) Z )
  -
  ( d^\ell e^m ) ( d Z )
  -
  ( d^\ell e^m ) [ Z, d ]
  +
  [ ( d^\ell e^m ) Z, d ].
  \]
Using induction we can write this as
  \[
  L( d^{\ell+1} e^m )
  =
  2
  L(d) L(d)^\ell L(e)^m
  -
  L(d)^\ell L(e)^m L(d)
  -
  [
  L(d)^\ell L(e)^m, \rho(d)
  ].
  \]
Corollary \ref{operatorcommutators} shows that the commutator is zero and that
the first and second terms combine. We now prove $L( c^k d^\ell e^m ) = L(c)^k
L(d)^\ell L(e)^m$ by induction on $k$. Lemma \ref{generalproduct} gives
  \[
  (c^{k+1} d^\ell e^m) Z
  =
  2
  c ( ( c^k d^\ell e^m ) Z )
  -
  ( c^k d^\ell e^m ) ( c Z )
  -
  ( c^k d^\ell e^m ) [ Z, c ]
  +
  [ ( c^k d^\ell e^m ) Z, c ].
  \]
By induction and Proposition \ref{abcde-brackets} we can write this
as
  \allowdisplaybreaks
  \begin{align*}
  L( c^{k+1} d^\ell e^m )
  &=
  2
  L(c) L(c)^k L(d)^\ell L(e)^m
  -
  L(c)^k L(d)^\ell L(e)^m L(c)
  \\
  &\qquad
  -
  L(c)^k L(d)^\ell L(e)^m \rho(c)
  +
  \rho(c) L(c)^k L(d)^\ell L(e)^m.
  \end{align*}
Corollary \ref{operatorcommutators} now gives
  \allowdisplaybreaks
  \begin{align*}
  L( c^{k+1} d^\ell e^m )
  &=
  2
  L(c)^{k+1} L(d)^\ell L(e)^m
  -
  L(c)^k \Big( L(c) L(d)^\ell - \ell L(d)^{\ell-1} L(e) \Big) L(e)^m
  \\
  &\qquad
  -
  L(c)^k \Big( \rho(c) L(d)^\ell + \ell L(d)^{\ell-1} L(e) \Big) L(e)^m
  +
  \rho(c) L(c)^k L(d)^\ell L(e)^m,
  \end{align*}
and cancelation completes the proof.  \end{proof}

\begin{proposition} \label{bcde-product}
We have
  \[
  L( b^j c^k d^\ell e^m )
  =
  \sum_{\alpha=0}^{\min(j,\ell)}
  \frac{\alpha!}{6^\alpha}
  \binom{j}{\alpha}
  \binom{\ell}{\alpha}
  D_a^{\alpha}
  L(b)^{j-\alpha}
  L(c)^k
  L(d)^{\ell-\alpha}
  L(e)^{m+\alpha}.
  \]
\end{proposition}

\begin{proof}
Induction on $j$; the basis $j = 0$ is Proposition \ref{cde-product}. We use
Corollary \ref{operatorcommutators} repeatedly. Lemma \ref{generalproduct}
gives
  \allowdisplaybreaks
  \begin{align*}
  &
  ( b^{j+1} c^k d^\ell e^m ) ( a^p b^q c^r d^s e^t )
  =
  \\
  &\qquad
  2
  b ( ( b^j c^k d^\ell e^m ) ( a^p b^q c^r d^s e^t ) )
  -
  ( b^j c^k d^\ell e^m) ( b ( a^p b^q c^r d^s e^t ) )
  \\
  &\qquad\qquad
  -
  ( b^j c^k d^\ell e^m ) [ a^p b^q c^r d^s e^t, b ]
  +
  [ ( b^j c^k d^\ell e^m )( a^p b^q c^r d^s e^t ), b ],
  \end{align*}
which we can write as $L( b^{j+1} c^k d^\ell e^m ) = A + B + C + D$ where
  \begin{alignat*}{2}
  A &=
  2
  L(b) L( b^j c^k d^\ell e^m ),
  &\qquad
  B &=
  -
  L( b^j c^k d^\ell e^m ) L(b),
  \\
  C &=
  -
  L( b^j c^k d^\ell e^m ) \rho(b),
  &\qquad
  D &=
  \rho(b) L( b^j c^k d^\ell e^m ).
  \end{alignat*}
Using induction and $[ D_a, L(b) ] = 0$ we get
  \[
  A =
  2
  \sum_{\alpha=0}^\ell
  \frac{\alpha!}{6^\alpha}
  \binom{j}{\alpha}
  \binom{\ell}{\alpha}
  D_a^{\alpha}
  L(b)^{j+1-\alpha}
  L(c)^k
  L(d)^{\ell-\alpha}
  L(e)^{m+\alpha}.
  \]
Using induction and $[D_a, L(d) ] = [ D_a, L(c) ] = 0$ we get $B = B' + B''$
where
  \allowdisplaybreaks
  \begin{align*}
  B'
  &=
  -
  \sum_{\alpha=0}^\ell
  \frac{\alpha!}{6^\alpha}
  \binom{j}{\alpha}
  \binom{\ell}{\alpha}
  D_a^{\alpha}
  L(b)^{j+1-\alpha}
  L(c)^k
  L(d)^{\ell-\alpha}
  L(e)^{m+\alpha},
  \\
  B''
  &=
  -
  \frac13
  \sum_{\alpha=0}^\ell
  \frac{\alpha!}{6^\alpha}
  \binom{j}{\alpha}
  (\ell{-}\alpha)
  \binom{\ell}{\alpha}
  D_a^{\alpha+1}
  L(b)^{j-\alpha}
  L(c)^k
  L(d)^{\ell-\alpha-1}
  L(e)^{m+\alpha+1}.
  \end{align*}
Using induction and $[ D_a, L(c) ] = 0$ we get $C = C' + C''$ where
  \allowdisplaybreaks
  \begin{align*}
  C'
  &=
  -
  \sum_{\alpha=0}^\ell
  \frac{\alpha!}{6^\alpha}
  \binom{j}{\alpha}
  \binom{\ell}{\alpha}
  D_a^{\alpha}
  \rho(b)
  L(b)^{j-\alpha}
  L(c)^k
  L(d)^{\ell-\alpha}
  L(e)^{m+\alpha},
  \\
  C'' &=
  \frac12
  \sum_{\alpha=0}^\ell
  \frac{\alpha!}{6^\alpha}
  \binom{j}{\alpha}
  (\ell{-}\alpha)
  \binom{\ell}{\alpha}
  D_a^{\alpha+1}
  L(b)^{j-\alpha}
  L(c)^k
  L(d)^{\ell-\alpha-1}
  L(e)^{m+\alpha+1}.
  \end{align*}
Using induction and $[ D_a, \rho(b) ] = 0$ we get
  \[
  D =
  \sum_{\alpha=0}^\ell
  \frac{\alpha!}{6^\alpha}
  \binom{j}{\alpha}
  \binom{\ell}{\alpha}
  D_a^{\alpha}
  \rho(b)
  L(b)^{j-\alpha}
  L(c)^k
  L(d)^{\ell-\alpha}
  L(e)^{m+\alpha}.
  \]
Terms $A$ and $B'$ combine to give
  \[
  A + B'
  =
  \sum_{\alpha=0}^\ell
  \frac{\alpha!}{6^\alpha}
  \binom{j}{\alpha}
  \binom{\ell}{\alpha}
  D_a^{\alpha}
  L(b)^{j+1-\alpha}
  L(c)^k
  L(d)^{\ell-\alpha}
  L(e)^{m+\alpha}.
  \]
Terms $B''$ and $C''$ combine to give
  \[
  B'' + C''
  =
  \sum_{\alpha=1}^{\ell+1}
  \frac{\alpha!}{6^\alpha}
  \binom{j}{\alpha{-}1}
  \binom{\ell}{\alpha}
  D_a^\alpha
  L(b)^{j+1-\alpha}
  L(c)^k
  L(d)^{\ell-\alpha}
  L(e)^{m+\alpha}.
  \]
Terms $C'$ and $D$ cancel, and then $A + B'$ and $B'' + C''$ combine using
Pascal's identity to give the result.  \end{proof}

Before we prove the formula for $L( a^i b^j c^k d^\ell e^m )$ we need a
straightening Lemma for moving $L(a) + \rho(a)$ through a product of operators.

\begin{definition} \label{Xmonomial}
Our standard order for a product of operators will be
  \[
  X
  =
  L(a)^{s}
  D_a^{t}
  L(b)^{u}
  D_b^{v}
  L(c)^{w}
  D_d^{x}
  L(d)^{y}
  L(e)^{z}.
  \]
Note that $D_x$ only appears for $x = a, b, d$. Furthermore, $L(a)$ precedes
$D_a$ and $L(b)$ precedes $D_b$, but $D_d$ precedes $L(d)$.
\end{definition}

\begin{lemma} \label{LR-straighten}
We have
  \allowdisplaybreaks
  \begin{align*}
  &
  L(a)^{s}
  D_a^{t}
  L(b)^{u}
  D_b^{v}
  L(c)^{w}
  D_d^{x}
  L(d)^{y}
  L(e)^{z}
  \big(L(a){+}\rho(a)\big)
  =
  \\
  &\qquad
  \big(L(a){+}\rho(a)\big)
  L(a)^{s}
  D_a^{t}
  L(b)^{u}
  D_b^{v}
  L(c)^{w}
  D_d^{x}
  L(d)^{y}
  L(e)^{z}
  \\
  &\qquad\qquad
  +
  t
  L(a)^{s}
  D_a^{t-1}
  L(b)^{u}
  D_b^{v}
  L(c)^{w}
  D_d^{x}
  L(d)^{y}
  L(e)^{z}
  \\
  &\qquad\qquad\qquad
  -
  \tfrac16
  u
  L(a)^{s}
  D_a^{t}
  L(b)^{u-1}
  D_b^{v}
  L(c)^{w}
  D_d^{x+1}
  L(d)^{y}
  L(e)^{z+1}
  \\
  &\qquad\qquad\qquad\qquad
  +
  \tfrac16
  y
  L(a)^{s}
  D_a^{t}
  L(b)^{u}
  D_b^{v+1}
  L(c)^{w}
  D_d^{x}
  L(d)^{y-1}
  L(e)^{z+1}.
  \end{align*}
Therefore
  \allowdisplaybreaks
  \begin{align*}
  &
  2 L(a) X - X L(a) - X \rho(a) + \rho(a) X =
  \\
  &\qquad
  L(a)^{s+1}
  D_a^{t}
  L(b)^{u}
  D_b^{v}
  L(c)^{w}
  D_d^{x}
  L(d)^{y}
  L(e)^{z}
  \\
  &\qquad\qquad
  -
  t
  L(a)^{s}
  D_a^{t-1}
  L(b)^{u}
  D_b^{v}
  L(c)^{w}
  D_d^{x}
  L(d)^{y}
  L(e)^{z}
  \\
  &\qquad\qquad\qquad
  +
  \tfrac16
  u
  L(a)^{s}
  D_a^{t}
  L(b)^{u-1}
  D_b^{v}
  L(c)^{w}
  D_d^{x+1}
  L(d)^{y}
  L(e)^{z+1}
  \\
  &\qquad\qquad\qquad\qquad
  -
  \tfrac16
  y
  L(a)^{s}
  D_a^{t}
  L(b)^{u}
  D_b^{v+1}
  L(c)^{w}
  D_d^{x}
  L(d)^{y-1}
  L(e)^{z+1}.
  \end{align*}
\end{lemma}

\begin{proof}
We write $R(a) = L(a) + \rho(a)$: the operator of right multiplication by $a$
in $U(\mathbb{M})$. Corollary \ref{operatorcommutators}, Lemma
\ref{DMcommutators} and Proposition \ref{abcde-brackets} give
  \begin{align*}
  &
  [ R(a), L(a) ] =
  [ R(a), D_b ] =
  [ R(a), L(c) ] =
  [ R(a), D_d ] =
  [ R(a), L(e) ] =
  0,
  \\
  &
  [ R(a), L(b) ] = \tfrac16 M_e D_d,
  \quad
  [ R(a), L(d) ] = - \tfrac16 M_e D_b,
  \quad
  [ R(a), D_a ] = -1.
  \end{align*}
These equations imply
  \allowdisplaybreaks
  \begin{align*}
  &
  \Big[
  L(a)^{s}
  D_a^{t}
  L(b)^{u}
  D_b^{v}
  L(c)^{w}
  D_d^{x}
  L(d)^{y}
  L(e)^{z},
  R(a)
  \Big]
  \\
  &\qquad
  =
  L(a)^{s}
  \Big[
  D_a^{t},
  R(a)
  \Big]
  L(b)^{u}
  D_b^{v}
  L(c)^{w}
  D_d^{x}
  L(d)^{y}
  L(e)^{z}
  \\
  &\qquad\qquad
  +
  L(a)^{s}
  D_a^{t}
  \Big[
  L(b)^{u},
  R(a)
  \Big]
  D_b^{v}
  L(c)^{w}
  D_d^{x}
  L(d)^{y}
  L(e)^{z}
  \\
  &\qquad\qquad\qquad
  +
  L(a)^{s}
  D_a^{t}
  L(b)^{u}
  D_b^{v}
  L(c)^{w}
  D_d^{x}
  \Big[
  L(d)^{y},
  R(a)
  \Big]
  L(e)^{z},
  \end{align*}
which gives the first equation. The second part follows easily.
\end{proof}

We use the multinomial coefficients
  \[
  \binom{n}{i_1,\hdots,i_k}
  =
  \frac{n!}{i_1! \cdots i_k! (n{-}i_1{-}\cdots{-}i_k)!},
  \]
with the convention that
  \[
  \binom{n}{i_1,\hdots,i_k} = 0
  \;
  \text{if either $i_j < 0$ for some $j$ or $\sum_j i_j > n$}.
  \]

\begin{proposition} \label{abcde-product}
We have
  \allowdisplaybreaks
  \begin{align*}
  &
  L( a^i b^j c^k d^\ell e^m )
  =
  \\
  &\quad
  \sum_{\alpha=0}^\ell
  \,
  \sum_{\beta=0}^{\min(i,\alpha)}
  \,
  \sum_{\gamma=0}^{\min(\ell-\alpha,\beta)}
  \,
  \sum_{\delta=0}^{\min(j-\alpha,\gamma)}
  (-1)^{\beta+\delta}
  \frac{\alpha!\beta!}{6^{\alpha+\gamma}}
  \binom{\alpha}{\beta{-}\gamma}
  \binom{i}{\beta}
  \binom{j}{\alpha,\delta}
  \binom{\ell}{\alpha,\gamma{-}\delta}
  \,\times
  \\
  &\quad\quad
  L(a)^{i-\beta}
  D_a^{\alpha-\beta+\gamma}
  L(b)^{j-\alpha-\delta}
  D_b^{\gamma-\delta}
  L(c)^k
  D_d^\delta
  L(d)^{\ell-\alpha-\gamma+\delta}
  L(e)^{m+\alpha+\gamma}.
  \end{align*}
\end{proposition}

\begin{proof}
Induction on $i$; the basis $i = 0$ is Proposition \ref{bcde-product}. For the
inductive step we use Lemma \ref{generalproduct} to get
  \allowdisplaybreaks
  \begin{align*}
  &
  ( a^{i+1} b^j c^k d^\ell e^m ) ( a^p b^q c^r d^s e^t )
  =
  \\
  &\qquad
  2
  a ( (a^i b^j c^k d^\ell e^m ) ( a^p b^q c^r d^s e^t ) )
  -
  (a^i b^j c^k d^\ell e^m) ( a ( a^p b^q c^r d^s e^t ) )
  \\
  &\qquad\qquad
  -
  ( a^i b^j c^k d^\ell e^m ) [ a^p b^q c^r d^s e^t, a ]
  +
  [ ( a^i b^j c^k d^\ell e^m )( a^p b^q c^r d^s e^t ), a ].
  \end{align*}
Therefore $L( a^{i+1} b^j c^k d^\ell e^m ) = A + B + C + D$ where
  \begin{alignat*}{2}
  A &=
  L(a) L( a^i b^j c^k d^\ell e^m ),
  &\quad
  B &=
  -
  L( a^i b^j c^k d^\ell e^m ) L(a),
  \\
  C &=
  -
  L( a^i b^j c^k d^\ell e^m ) \rho(a),
  &\quad
  D &=
  \rho(a) L( a^i b^j c^k d^\ell e^m ).
  \end{alignat*}
We apply the second part of Lemma \ref{LR-straighten} to the
monomial
  \[
  X =
  L(a)^{i-\beta}
  D_a^{\alpha-\beta+\gamma}
  L(b)^{j-\alpha-\delta}
  D_b^{\gamma-\delta}
  L(c)^k
  D_d^\delta
  L(d)^{\ell-\alpha-\gamma+\delta}
  L(e)^{m+\alpha+\gamma}.
  \]
Therefore
  \allowdisplaybreaks
  \begin{align*}
  A
  &=
  \sum_{\alpha=0}^\ell
  \sum_{\beta=0}^i
  \sum_{\gamma=0}^\beta
  \sum_{\delta=0}^\gamma
  (-1)^{\beta+\delta}
  \frac{\alpha!\beta!}{6^{\alpha+\gamma}}
  P
  Q
  R
  S
  \;\times
  \\
  &
  L(a)^{i-\beta+1}
  D_a^{\alpha-\beta+\gamma}
  L(b)^{j-\alpha-\delta}
  D_b^{\gamma-\delta}
  L(c)^{k}
  D_d^{\delta}
  L(d)^{\ell-\alpha-\gamma+\delta}
  L(e)^{m+\alpha+\gamma},
  \\
  B
  &=
  \sum_{\alpha=0}^\ell
  \sum_{\beta=0}^i
  \sum_{\gamma=0}^\beta
  \sum_{\delta=0}^\gamma
  (-1)^{\beta+\delta+1}
  \frac{\alpha!\beta!}{6^{\alpha+\gamma}}
  (\alpha{-}\beta{+}\gamma)
  P
  Q
  R
  S
  \;\times
  \\
  &
  L(a)^{i-\beta}
  D_a^{\alpha-\beta+\gamma-1}
  L(b)^{j-\alpha-\delta}
  D_b^{\gamma-\delta}
  L(c)^{k}
  D_d^{\delta}
  L(d)^{\ell-\alpha-\gamma+\delta}
  L(e)^{m+\alpha+\gamma},
  \\
  C
  &=
  \sum_{\alpha=0}^\ell
  \sum_{\beta=0}^i
  \sum_{\gamma=0}^\beta
  \sum_{\delta=0}^\gamma
  (-1)^{\beta+\delta}
  \frac{\alpha!\beta!}{6^{\alpha+\gamma+1}}
  P
  Q
  (j{-}\alpha{-}\delta)
  R
  S
  \;\times
  \\
  &
  L(a)^{i-\beta}
  D_a^{\alpha-\beta+\gamma}
  L(b)^{j-\alpha-\delta-1}
  D_b^{\gamma-\delta}
  L(c)^{k}
  D_d^{\delta+1}
  L(d)^{\ell-\alpha-\gamma+\delta}
  L(e)^{m+\alpha+\gamma+1},
  \\
  D
  &=
  \sum_{\alpha=0}^\ell
  \sum_{\beta=0}^i
  \sum_{\gamma=0}^\beta
  \sum_{\delta=0}^\gamma
  (-1)^{\beta+\delta+1}
  \frac{\alpha!\beta!}{6^{\alpha+\gamma+1}}
  P
  Q
  R
  (\ell{-}\alpha{-}\gamma{+}\delta)
  S
  \;\times
  \\
  &
  L(a)^{i-\beta}
  D_a^{\alpha-\beta+\gamma}
  L(b)^{j-\alpha-\delta}
  D_b^{\gamma-\delta+1}
  L(c)^{k}
  D_d^{\delta}
  L(d)^{\ell-\alpha-\gamma+\delta-1}
  L(e)^{m+\alpha+\gamma+1},
  \end{align*}
where
  \[
  P = \binom{\alpha}{\beta{-}\gamma},
  \qquad
  Q = \binom{i}{\beta},
  \qquad
  R = \binom{j}{\alpha,\delta},
  \qquad
  S = \binom{\ell}{\alpha,\gamma{-}\delta}.
  \]
We make the following substitutions in the summation indices: in $B$ we replace
$\beta$ by $\beta{-}1$; in $C$ we replace $\beta$ by $\beta{-}1$, $\gamma$ by
$\gamma{-}1$, and $\delta$ by $\delta{-}1$; in $D$ we replace $\beta$ by
$\beta{-}1$, and $\gamma$ by $\gamma{-}1$. Using our convention on multinomial
coefficients, we can write all four sums with the notation
  \allowdisplaybreaks
  \begin{align*}
  {\sum}^4
  &=
  \sum_{\alpha=0}^\ell
  \sum_{\beta=0}^{i+1}
  \sum_{\gamma=0}^\beta
  \sum_{\delta=0}^\gamma
  (-1)^{\beta+\delta}
  \frac{\alpha!}{6^{\alpha+\gamma}},
  \\
  Y
  &=
  L(a)^{i-\beta+1}
  D_a^{\alpha-\beta+\gamma}
  L(b)^{j-\alpha-\delta}
  D_b^{\gamma-\delta}
  L(c)^{k}
  D_d^{\delta}
  L(d)^{\ell-\alpha-\gamma+\delta}
  L(e)^{m+\alpha+\gamma}.
  \end{align*}
We obtain
  \allowdisplaybreaks
  \begin{align*}
  A
  &=
  {\sum}^4
  \big[
  \beta!
  P
  Q
  R
  S
  \big]
  Y,
  \\
  B
  &=
  {\sum}^4
  \bigg[
  (\beta{-}1)!
  (\alpha{-}\beta{+}\gamma{+}1)
  \binom{\alpha}{\beta{-}\gamma{-}1}
  \binom{i}{\beta{-}1}
  R
  S
  \bigg]
  Y,
  \\
  C
  &=
  {\sum}^4
  \bigg[
  (\beta{-}1)!
  P
  \binom{i}{\beta{-}1}
  (j{-}\alpha{-}\delta{+}1)
  \binom{j}{\alpha,\delta{-}1}
  S
  \bigg]
  Y,
  \\
  D
  &=
  {\sum}^4
  \bigg[
  (\beta{-}1)!
  P
  \binom{i}{\beta{-}1}
  R
  (\ell{-}\alpha{-}\gamma{+}\delta{+}1)
  \binom{\ell}{\alpha,\gamma{-}\delta{-}1}
  \bigg]
  Y.
  \end{align*}
The sum of these four terms is
  \[
  {\sum}^4
  \beta!
  \bigg[
  \frac{i{+}1{-}\beta}{i{+}1}
  {+}
  \frac{\beta{-}\gamma}{i{+}1}
  {+}
  \frac{\delta}{i{+}1}
  {+}
  \frac{\gamma{-}\delta}{i{+}1}
  \bigg]
  P \binom{i{+}1}{\beta} R S
  \,
  Y
  =
  {\sum}^4
  \bigg[
  \beta!
  P \binom{i{+}1}{\beta} R S
  \bigg]
  Y,
  \]
as required.
\end{proof}

\begin{lemma} \label{bd-expansions}
The powers of $L(b)$ and $L(d)$ are
  \allowdisplaybreaks
  \begin{align*}
  L(b)^u
  &=
  \sum_{\epsilon=0}^u
  \sum_{\zeta=0}^{u-\epsilon}
  (-1)^\zeta
  \frac{1}{3^{u-\epsilon-\zeta}}
  \binom{u}{\epsilon,\zeta}
  M_b^\epsilon
  M_c^\zeta
  M_e^{u-\epsilon-\zeta}
  D_a^{u-\epsilon}
  D_d^{u-\epsilon-\zeta}
  ,
  \\
  L(d)^y
  &=
  \sum_{\eta=0}^y
  \sum_{\theta=0}^{y-\eta}
  (-1)^{y-\eta}
  \frac{1}{3^{y-\eta-\theta}}
  \binom{y}{\eta,\theta}
  M_d^\eta
  M_e^{y-\eta}
  D_a^{y-\eta-\theta}
  D_b^{y-\eta-\theta}
  D_c^\theta.
  \end{align*}
\end{lemma}

\begin{proof}
We use the trinomial theorem, since the terms in $L(b)$ and $L(d)$ commute.

\end{proof}

\begin{lemma} \label{finalstraightening}
The expansion of the monomial $X$ of Definition \ref{Xmonomial} is
  \allowdisplaybreaks
  \begin{align*}
  &
  L(a)^{s}
  D_a^{t}
  L(b)^{u}
  D_b^{v}
  L(c)^{w}
  D_d^{x}
  L(d)^{y}
  L(e)^{z}
  =
  \\
  &
  \sum_{\epsilon=0}^u
  \sum_{\zeta=0}^{u-\epsilon}
  \sum_{\eta=0}^y
  \sum_{\theta=0}^{y-\eta}
  \,
  \sum_{\lambda=0}^{\min(u-\epsilon-\zeta+x,\eta)}
  \frac{(-1)^{\zeta+y-\eta}\lambda!}{3^{u-\epsilon-\zeta+y-\eta-\theta}}
  \binom{u}{\epsilon,\zeta}
  \binom{y}{\eta,\theta}
  \binom{u{-}\epsilon{-}\zeta{+}x}{\lambda}
  \binom{\eta}{\lambda}
  \\
  &
  M_a^s
  M_b^\epsilon
  M_c^{\zeta+w}
  M_d^{\eta-\lambda}
  M_e^{u-\epsilon-\zeta+y-\eta+z}
  D_a^{t+u-\epsilon+y-\eta-\theta}
  D_b^{v+y-\eta-\theta}
  D_c^\theta
  D_d^{u+x-\epsilon-\zeta-\lambda}.
  \end{align*}
\end{lemma}

\begin{proof}
Proposition \ref{abcde-brackets} and Lemma \ref{bd-expansions} give
  \allowdisplaybreaks
  \begin{align*}
  &
  L(a)^{s}
  D_a^{t}
  L(b)^{u}
  D_b^{v}
  L(c)^{w}
  D_d^{x}
  L(d)^{y}
  L(e)^{z}
  =
  \\
  &
  \sum_{\epsilon=0}^u
  \sum_{\zeta=0}^{u-\epsilon}
  \sum_{\eta=0}^y
  \sum_{\theta=0}^{y-\eta}
  \frac{(-1)^{\zeta+y-\eta}}{3^{u-\epsilon-\zeta+y-\eta-\theta}}
  \binom{u}{\epsilon,\zeta}
  \binom{y}{\eta,\theta}
  \times
  \\
  &
  M_a^s
  D_a^t
  M_b^\epsilon
  M_c^\zeta
  M_e^{u-\epsilon-\zeta}
  D_a^{u-\epsilon}
  D_d^{u-\epsilon-\zeta}
  D_b^v
  M_c^w
  D_d^x
  M_d^\eta
  M_e^{y-\eta}
  D_a^{y-\eta-\theta}
  D_b^{y-\eta-\theta}
  D_c^\theta
  M_e^z.
  \end{align*}
Using Lemma \ref{DMcommutators} we can move $M_x$ to the left and $D_x$ to the
right, and collect the remaining noncommuting factors on the right:
  \[
  M_a^s
  M_b^\epsilon
  M_c^{\zeta+w}
  M_e^{u-\epsilon-\zeta+y-\eta+z}
  D_a^{t+u-\epsilon+y-\eta-\theta}
  D_b^{v+y-\eta-\theta}
  D_c^\theta
  \big(
  D_d^{u-\epsilon-\zeta+x}
  M_d^\eta
  \big)
  \]
To complete the proof we use the commutation formula
  \[
  D_d^m M_d^n
  =
  \sum_{i=0}^{\min(m,n)}
  i!
  \binom{m}{i}
  \binom{n}{i}
  M_d^{n-i} D_d^{m-i}.
  \]
(Compare Proposition \ref{nilpotentLiealgebra}.)  \end{proof}

\begin{proposition} \label{finaloperatorexpansion}
We have
  \allowdisplaybreaks
  \begin{align*}
  &
  L( a^i b^j c^k d^\ell e^m )
  =
  \\
  &\qquad
  \sum_{\alpha=0}^\ell \;
  \sum_{\beta=0}^i \;
  \sum_{\gamma=0}^\beta \;
  \sum_{\delta=0}^\gamma \;
  \sum_{\epsilon=0}^{j-\alpha-\delta} \;
  \sum_{\zeta=0}^{j-\alpha-\delta-\epsilon} \;
  \sum_{\eta=0}^{\ell-\alpha-(\gamma-\delta)} \;
  \sum_{\theta=0}^{\ell-\alpha-(\gamma-\delta)-\eta} \;
  \sum_{\lambda=0}^\eta
  \\
  &\qquad\qquad
  (-1)^{\beta+\zeta+\ell-\alpha-\gamma-\eta}
  \frac{\alpha!\beta!\lambda!}
  {2^{\alpha+\gamma}3^{j-\epsilon-\zeta+\ell-\alpha-\eta-\theta}}
  \,\times
  \\
  &\qquad\qquad\qquad
  \binom{\alpha}{\beta{-}\gamma}
  \binom{i}{\beta}
  \binom{j}{\alpha,\delta,\epsilon,\zeta}
  \binom{j{-}\alpha{-}\epsilon{-}\zeta}{\lambda}
  \binom{\ell}{\alpha,\gamma{-}\delta,\eta,\theta}
  \binom{\eta}{\lambda}
  \,\times
  \\
  &\qquad\qquad\qquad\qquad
  M_a^{i-\beta}
  M_b^\epsilon
  M_c^{\zeta+k}
  M_d^{\eta-\lambda}
  M_e^{j-\alpha-\epsilon-\zeta+\ell-\eta+m}
  \,\times
  \\
  &\qquad\qquad\qquad\qquad\qquad
  D_a^{j-\beta-\epsilon+\ell-\alpha-\eta-\theta}
  D_b^{\ell-\alpha-\eta-\theta}
  D_c^\theta
  D_d^{j-\alpha-\epsilon-\zeta-\lambda}.
  \end{align*}
\end{proposition}

\begin{proof}
In Proposition \ref{abcde-product} take $s = i{-}\beta$, $t =
\alpha{-}\beta{+}\gamma$, $u = j{-}\alpha{-}\delta$, $v = \gamma{-}\delta$, $w
= k$, $x = \delta$, $y = \ell{-}\alpha{-}\gamma{+}\delta$, $z =
m{+}\alpha{+}\gamma$, and apply Lemma \ref{finalstraightening}.
\end{proof}

In the next result we use the notation
  \[
  \begin{bmatrix} n \\ 0 \end{bmatrix}
  =
  1,
  \;
  \begin{bmatrix} n \\ k \end{bmatrix}
  =
  n(n{-}1)\cdots(n{-}k{+}1),
  \;
  \text{so that}
  \;
  D_x^k ( x^n ) =
  \begin{bmatrix} n \\ k \end{bmatrix}
  x^{n-k}.
  \]

\begin{theorem}[Universal structure constants] \label{maintheorem}
In $U(\mathbb{M})$ we have
  \allowdisplaybreaks
  \begin{align*}
  &
  ( a^i b^j c^k d^\ell e^m ) ( a^p b^q c^r d^s e^t )
  =
  \\
  &
  \sum_{\alpha=0}^\ell
  \sum_{\beta=0}^i
  \sum_{\gamma=0}^\beta
  \sum_{\delta=0}^\gamma
  \sum_{\epsilon=0}^{j-\alpha-\delta}
  \sum_{\zeta=0}^{j-\alpha-\delta-\epsilon}
  \sum_{\eta=0}^{\ell-\alpha-(\gamma-\delta)}
  \sum_{\theta=0}^{\ell-\alpha-(\gamma-\delta)-\eta}
  \sum_{\lambda=0}^\eta
  \\
  &
  (-1)^{\beta+\zeta+\ell-\alpha-\gamma-\eta}
  \frac{\alpha!\beta!\lambda!}{2^{\alpha+\gamma}3^{j-\epsilon-\zeta+\ell-\alpha-\eta-\theta}}
  \,\times
  \\
  &
  \binom{\alpha}{\beta{-}\gamma}
  \binom{i}{\beta}
  \binom{j}{\alpha,\delta,\epsilon,\zeta}
  \binom{j{-}\alpha{-}\epsilon{-}\zeta}{\lambda}
  \binom{\ell}{\alpha,\gamma{-}\delta,\eta,\theta}
  \binom{\eta}{\lambda}
  \,\times
  \\
  &
  \begin{bmatrix}
  p \\ j{-}\beta{-}\epsilon{+}\ell{-}\alpha{-}\eta{-}\theta
  \end{bmatrix}
  \begin{bmatrix}
  q \\ \ell{-}\alpha{-}\eta{-}\theta
  \end{bmatrix}
  \begin{bmatrix}
  r \\ \theta
  \end{bmatrix}
  \begin{bmatrix}
  s \\ j{-}\alpha{-}\epsilon{-}\zeta{-}\lambda
  \end{bmatrix}
  \,\times
  \\
  &
  a^{i{+}p{-}j{+}\epsilon{-}\ell{+}\alpha{+}\eta{+}\theta}
  b^{\epsilon{+}q{-}\ell{+}\alpha{+}\eta{+}\theta}
  c^{\zeta{+}k{+}r{-}\theta}
  d^{\eta{+}s{+}\alpha{+}\epsilon{+}\zeta{-}j}
  e^{j{-}\alpha{-}\epsilon{-}\zeta{+}\ell{-}\eta{+}m{+}t}.
  \end{align*}
\end{theorem}


\section{The universal alternative enveloping algebra}

The algebra $U(\mathbb{M})$ is not power-associative, since
  \[
  ( abd, abd, abd )
  =
  \tfrac16 a b c d^2 e
  -
  \tfrac16 a b d e^2
  -
  \tfrac16 c^2 d^2 e
  +
  \tfrac{11}{36} c d e^2
  -
  \tfrac{1}{12} e^3.
  \]
In this section we construct the maximal alternative quotient of
$U(\mathbb{M})$.

\begin{definition}
The \textbf{alternator ideal} of a nonassociative algebra $A$ is
  \[
  I(A) = \langle \, (x,x,y), \, (y,x,x) \, | \, x, y \in A \, \rangle
  \]
If $M$ is a Malcev algebra then we write $I(M) = I(U(M))$. The
\textbf{universal alternative enveloping algebra} of $M$ is $A(M) = U(M)/I(M)$.
\end{definition}

\begin{remark}
The speciality problem for Malcev algebra is equivalent to the question of the
injectivity of the natural mapping $M \to A(M)$.
\end{remark}

In $U(\mathbb{M})$ we have the alternators
  \[
  (ab,ab,d) = -\tfrac16 ce,
  \qquad
  (bd,bd,a^2) = \tfrac{1}{18} e^2.
  \]

\begin{definition} \label{idealJdefinition}
We write $I = I(\mathbb{M})$ and $J = \mathrm{ideal}\langle ce, e^2 \rangle$;
it is clear that $J \subseteq I$.
\end{definition}

\begin{lemma} \label{refereelemma}
A basis of the ideal $J$ consists of the set of monomials
  \[
  \{ a^i b^j c^k d^\ell e^m \,|\, m \ge 2 \}
  \cup
  \{ a^i b^j c^k d^\ell e \,|\, k \ge 1 \}.
  \]
\end{lemma}

\begin{proof}
Linear independence is clear. Every monomial $a^i b^j c^k d^\ell e^m$ ($m \ge
2$) belongs to $J$. Since $c d^\ell e = \ell d^{\ell-1} e^2 + d^\ell ce$, every
monomial $a^i b^j c^k d^\ell e$ ($k \ge 1$) belongs to $J$. Since the
generators of $J$ belong to this set, it suffices to show that the span of
these monomials is an ideal in $U(\mathbb{M})$. This follows from Theorem
\ref{maintheorem}: the product of $a^i b^j c^k d^\ell e^m$ and $a^p b^q c^r d^s
e^t$ is a linear combination of monomials in which the exponents of $e$ satisfy
  \begin{align*}
  &
  j - \alpha - \epsilon - \zeta + \ell - \eta + m + t
  \ge
  \delta + \ell - \eta + m + t
  \ge
  \alpha + \gamma + m + t
  \ge
  m + t.
  \end{align*}
This exponent is 0 if and only if $m = t = 0$. For $(m,t) \ne (0,0)$, this
exponent is 1 if and only if $(m,t) \in \{ (1,0), (0,1) \}$, $\delta = j -
\alpha - \epsilon - \zeta$, and $\delta + \ell - \eta = \alpha + \gamma$. This
forces $\theta = 0$, and then the exponent of $c$ is $k + r + \zeta \ge k + r$.
It follows that the product of an arbitrary monomial with a monomial in the set
is a linear combination of monomials in the set.
\end{proof}

We now determine the structure constants for the quotient algebra
$U(\mathbb{M})/J$. We will show that $U(\mathbb{M})/J$ is an alternative
algebra; it will then follow that $I = J$ and that $A(\mathbb{M}) =
U(\mathbb{M})/J$. A spanning set for $U(\mathbb{M})/J$ (in fact a basis, by
Lemma \ref{refereelemma}) consists of the cosets of the monomials $m = a^i b^j
d^\ell e$ (type 1: the exponent of $e$ is 1, and so the exponent of $c$ is 0)
and $m = a^i b^j c^k d^\ell$ (type 2: the exponent of $e$ is 0). For type 1,
since $\{ c, d, e \}$ spans a Lie subalgebra of $\mathbb{M}$, we have $[ c,
d^\ell ] = \ell d^{\ell-1} e$ and so $c d^\ell e = \ell d^{\ell-1} e^2 + d^\ell
c e$; thus $c$ cannot occur. In the next result we write $m$ for the coset $m +
J$.

\begin{theorem}[Alternative structure constants] \label{alternativequotient}
In $U(\mathbb{M})/J$ we have
  \allowdisplaybreaks
  \begin{align}
  &
  ( a^i b^j d^\ell e ) ( a^p b^q d^s e )
  =
  0,
  \label{1times1}
  \\
  &
  ( a^i b^j c^k d^\ell ) ( a^p b^q d^s e )
  =
  \delta_{k0}
  a^{i+p} b^{j+q} d^{\ell+s} e,
  \label{1times2}
  \\
  &
  ( a^i b^j d^\ell e ) ( a^p b^q c^r d^s )
  =
  \delta_{r0}
  a^{i+p} b^{j+q} d^{\ell+s} e,
  \label{2times1}
  \\
  &
  ( a^i b^j c^k d^\ell ) ( a^p b^q c^r d^s )
  =
  \label{2times2}
  \sum_{\mu=0}^j
  (-1)^\mu
  \mu!
  \binom{j}{\mu}
  \binom{p}{\mu}
  a^{i+p-\mu}
  b^{j+q-\mu}
  c^{k+r+\mu}
  d^{\ell+s}
  \\
  &\qquad
  +
  \delta_{k0} \delta_{r0}
  \left(
    \tfrac16 i j s
  {-} \tfrac16 i \ell q
  {+} \tfrac12 j \ell p
  {+} \tfrac13 j p s
  {-} \tfrac13 \ell p q
  \right)
  a^{i+p-1} b^{j+q-1} d^{\ell+s-1} e
  \notag
  \\
  &\qquad\qquad
  - \delta_{k0} \delta_{r1} \ell a^{i+p} b^{j+q} d^{\ell+s-1} e.
  \notag
  \end{align}
\end{theorem}

\begin{proof}
We need only the terms on the right side of Theorem \ref{maintheorem} which are
nonzero modulo $J$:
  \begin{alignat*}{2}
  &\text{either}
  &\qquad
  &
  j{-}\alpha{-}\epsilon{-}\zeta{+}\ell{-}\eta{+}m{+}t = 0,
  \\
  &\text{or}
  &\qquad
  &
  j{-}\alpha{-}\epsilon{-}\zeta{+}\ell{-}\eta{+}m{+}t = 1,
  \quad
  \zeta{+}k{+}r{-}\theta = 0.
  \end{alignat*}
We write the exponent of $e$ as the sum of three nonnegative terms:
  \begin{equation}
  (j{-}\alpha{-}\epsilon{-}\zeta{-}\delta)
  +
  (\ell{-}\eta{+}\delta{-}\alpha{-}\gamma)
  +
  (\alpha{+}\gamma{+}m{+}t).
  \label{nonnegative}
  \end{equation}
Equation (\ref{1times1}): We have $m = t = 1$: on the right side of Theorem
\ref{maintheorem}, the exponent of $e$ in every term is $\ge 2$, so every term
becomes zero in $U(\mathbb{M})/J$.

\noindent Equation (\ref{1times2}): We have $m = 0$, $r = 0$ and $t = 1$: the
exponent of $e$ must be 1 and hence the exponent of $c$ must be 0. Therefore
  \begin{align*}
  &
  j{-}\alpha{-}\epsilon{-}\zeta{-}\delta = 0,
  \quad
  \ell{-}\eta{+}\delta{-}\alpha{-}\gamma = 0,
  \quad
  \alpha = 0,
  \quad
  \gamma = 0,
  \quad
  \zeta{+}k{-}\theta = 0.
  \end{align*}
The factor $\binom{\alpha}{\beta-\gamma}$ will be zero unless $\beta = 0$. From
the sum on $\delta$ we get $\delta = 0$. We now get $j{-}\epsilon{-}\zeta = 0$
and $\ell{-}\eta = 0$, so $\zeta = j{-}\epsilon$ and $\eta = \ell$. From the
sum on $\theta$ we get $\theta = 0$. Hence $\zeta{+}k = 0$ and so $\zeta = k =
0$; then $\epsilon = j$. The factor
$\binom{j{-}\alpha{-}\epsilon{-}\zeta}{\lambda}$ will be zero unless $\lambda =
0$. Only one term remains.

\noindent Equation (\ref{2times1}): Similar to Equation (\ref{1times2}).

\noindent Equation (\ref{2times2}): In Theorem \ref{maintheorem} we have $m =
0$ and $t = 0$. We write the result as $T_0 + T_1$, collecting terms with the
same exponent of $e$. For $T_0$, all three terms in (\ref{nonnegative}) must be
0. Since $\alpha{+}\gamma = 0$ we get $\alpha = \gamma = 0$, and hence $\beta =
\delta = 0$. Then $j{-}\epsilon{-}\zeta = 0$ and $\ell{-}\eta = 0$, which imply
$\zeta = j{-}\epsilon$ and $\eta = \ell$. Then $\theta = 0$, and the factor
$\binom{j{-}\alpha{-}\epsilon{-}\zeta}{\lambda}$ will be zero unless $\lambda =
0$. The remaining terms correspond to Proposition \ref{nilpotentLiealgebra}
with the new factor $d^{\ell+s}$:
  \allowdisplaybreaks
  \begin{align*}
  T_0
  &=
  \sum_{\epsilon=0}^j
  (-1)^{j-\epsilon}
  \binom{j}{\epsilon}
  \begin{bmatrix}
  p \\ j{-}\epsilon
  \end{bmatrix}
  a^{i-(j-\epsilon)+p}
  b^{\epsilon+q}
  c^{j-\epsilon+k+r}
  d^{\ell+s}
  \\
  &=
  \sum_{\mu=0}^j
  (-1)^\mu
  \mu!
  \binom{j}{\mu}
  \binom{p}{\mu}
  a^{i+p-\mu}
  b^{j+q-\mu}
  c^{k+r+\mu}
  d^{\ell+s}.
  \end{align*}
For $T_1$, we have
  \[
  (j{-}\alpha{-}\epsilon{-}\zeta{-}\delta)
  +
  (\ell{-}\eta{+}\delta{-}\alpha{-}\gamma)
  +
  (\alpha{+}\gamma)
  =
  1,
  \quad
  \zeta{+}k{+}r{-}\theta = 0.
  \]
The four cases
  \begin{align*}
  &(1) \quad
  j{-}\alpha{-}\epsilon{-}\zeta{-}\delta = 1, \quad
  \ell{-}\eta{+}\delta{-}\alpha{-}\gamma = 0, \quad
  \alpha = 0, \quad
  \gamma = 0
  \\
  &(2) \quad
  j{-}\alpha{-}\epsilon{-}\zeta{-}\delta = 0, \quad
  \ell{-}\eta{+}\delta{-}\alpha{-}\gamma = 1, \quad
  \alpha = 0, \quad
  \gamma = 0
  \\
  &(3) \quad
  j{-}\alpha{-}\epsilon{-}\zeta{-}\delta = 0, \quad
  \ell{-}\eta{+}\delta{-}\alpha{-}\gamma = 0, \quad
  \alpha = 1, \quad
  \gamma = 0
  \\
  &(4) \quad
  j{-}\alpha{-}\epsilon{-}\zeta{-}\delta = 0, \quad
  \ell{-}\eta{+}\delta{-}\alpha{-}\gamma = 0, \quad
  \alpha = 0, \quad
  \gamma = 1
  \end{align*}
produce respectively 1, 2, 1, 2 terms, giving
  \begin{align*}
  T_1
  &=
  \delta_{k0} \delta_{r0}
  \left(
    \tfrac16 i j s
  - \tfrac16 i \ell q
  + \tfrac12 j \ell p
  + \tfrac13 j p s
  - \tfrac13 \ell p q
  \right)
  a^{i+p-1} b^{j+q-1} d^{\ell+s-1} e
  \\
  &\quad
  - \delta_{k0} \delta_{r1} \ell a^{i+p} b^{j+q} d^{\ell+s-1} e.
  \end{align*}
This completes the proof.
\end{proof}

\begin{corollary} \label{zeroassociators}
We have $(m,m',m'') = (m',m,m'') = (m',m'',m) = 0$ in $U(\mathbb{M})/J$, where
$m$ has type 1 and $m'$, $m''$ are arbitrary.
\end{corollary}

\begin{proof}
This follows easily from Theorem \ref{alternativequotient}.
\end{proof}

\begin{corollary} \label{type2associator}
The associator of type 2 monomials is
  \begin{align*}
  &
  ( a^i b^j c^k d^\ell, a^p b^q c^r d^s, a^v b^w c^x d^y )
  =
  \\
  &
  \delta_{k0} \delta_{r0} \delta_{x0}
  \tfrac16
  \big(
  i q y {-} i s w {-} j p y {+} j s v {+} \ell p w {-} \ell q v
  \big)
  a^{i+p+v-1} b^{j+q+w-1} d^{\ell+s+y-1} e.
  \end{align*}
\end{corollary}

\begin{proof}
We write $m_1 = a^i b^j c^k d^\ell$, $m_2 = a^p b^q c^r d^s$, $m_3 = a^v b^w
c^x d^y$. We use the notation $N( a^i b^j c^k, a^p b^q c^r )$ for the
associative multiplication in the enveloping algebra of the nilpotent Lie
subalgebra of $\mathbb{M}$ with basis $\{a,b,c\}$ (compare Proposition
\ref{nilpotentLiealgebra}); we extend $N$ in the obvious way to linear
combinations of monomials. Theorem \ref{alternativequotient} shows that
  \[
  ( a^i b^j c^k d^\ell ) ( a^p b^q c^r d^s )
  =
  N( a^i b^j c^k, a^p b^q c^r ) d^{\ell+s}
  +
  T_1(i,j,k,\ell,p,q,r,s),
  \]
where $T_1$ denotes the terms involving $e$. We define
  \[
  C(i,j,\ell,p,q,s)
  =
  \tfrac16 i j s
  -
  \tfrac16 i \ell q
  +
  \tfrac12 j \ell p
  +
  \tfrac13 j p s
  -
  \tfrac13 \ell p q.
  \]
We obtain
  \allowdisplaybreaks
  \begin{align*}
  ( m_1 m_2 ) m_3
  &=
  \big( ( a^i b^j c^k d^\ell ) ( a^p b^q c^r d^s ) \big)
  ( a^v b^w c^x d^y )
  \\
  &=
  \big( N( a^i b^j c^k, a^p b^q c^r ) d^{\ell+s} \big) ( a^v b^w c^x d^y )
  +
  T_1(i,j,k,\ell,p,q,r,s) ( a^v b^w c^x d^y )
  \\
  &=
  N( N( a^i b^j c^k, a^p b^q c^r ), a^v b^w c^x ) d^{\ell+s+y}
  \\
  &\quad
  +
  \delta_{k0} \delta_{r0} \delta_{x0}
  \big[
  C(i{+}p,j{+}q,\ell{+}s,v,w,y)
  +
  C(i,j,\ell,p,q,s)
  \big]
  \times
  \\
  &\quad\quad\quad
  a^{i+p+v-1} b^{j+q+w-1} d^{\ell+s+y-1} e
  \\
  &\quad
  -
  \delta_{k0} \delta_{r0} \delta_{x1}
  (\ell{+}s)
  a^{i+p+v} b^{j+q+w} d^{\ell+s+y-1} e
  \\
  &\quad
  -
  \delta_{k0} \delta_{r1} \delta_{x0}
  \ell
  a^{i+p+v} b^{j+q+w} d^{\ell+s+y-1} e,
  \\
  m_1 ( m_2 m_3 )
  &=
  ( a^i b^j c^k d^\ell )
  \big( ( a^p b^q c^r d^s ) ( a^v b^w c^x d^y ) \big)
  \\
  &=
  ( a^i b^j c^k d^\ell )
  \big( N( a^p b^q c^r, a^v b^w c^x ) d^{s+y} \big)
  +
  ( a^i b^j c^k d^\ell )
  T_1(p,q,r,s,v,w,x,y)
  \\
  &=
  N( a^i b^j c^k, N( a^p b^q c^r, a^v b^w c^x ) ) d^{\ell+s+y}
  \\
  &\quad
  +
  \delta_{k0} \delta_{r0} \delta_{x0}
  \big[
  C(i,j,\ell,p{+}v,q{+}w,s{+}y)
  +
  \ell q v
  +
  C(p,q,s,v,w,y)
  \big]
  \times
  \\
  &\quad\quad\quad
  a^{i+p+v-1} b^{j+q+w-1} d^{\ell+s+y-1} e
  \\
  &\quad
  -
  \delta_{k0} \delta_{r1} \delta_{x0}
  \ell
  a^{i+p+v} b^{j+q+w} d^{\ell+s+y-1} e
  \\
  &\quad
  -
  \delta_{k0} \delta_{r0} \delta_{x1}
  (\ell{+}s)
  a^{i+p+v} b^{j+q+w} d^{\ell+s+y-1} e.
  \end{align*}
The associator $( m_1, m_2, m_3 )$ is therefore
  \allowdisplaybreaks
  \begin{align*}
  \delta_{k0} \delta_{r0} \delta_{x0}
  &\big[
  C(i{+}p,j{+}q,\ell{+}s,v,w,y)
  +
  C(i,j,\ell,p,q,s)
  \\
  &\qquad
  -
  C(i,j,\ell,p{+}v,q{+}w,s{+}y)
  -
  \ell q v
  -
  C(p,q,s,v,w,y)
  \big]
  \times
  \\
  &
  a^{i+p+v-1} b^{j+q+w-1} d^{\ell+s+y-1} e.
  \end{align*}
The expression in square brackets simplifies as required.
\end{proof}

\begin{corollary}
The algebra $U(\mathbb{M})/J$ is alternative.
\end{corollary}

\begin{proof}
Corollaries \ref{zeroassociators} and \ref{type2associator} show that the
associator alternates.
\end{proof}

\begin{corollary}
The alternator ideal $I(\mathbb{M})$ is generated by $ce$ and $e^2$.
\end{corollary}

\begin{corollary}
The universal alternative enveloping algebra $A(\mathbb{M})$ is isomorphic to
the algebra with basis $\{ \, a^i b^j c^k d^\ell, \, a^i b^j d^\ell e \, | \,
i, j, k, \ell \ge 0 \, \}$ and structure constants of Proposition
\ref{alternativequotient}.
\end{corollary}

\begin{proof}
This follows from Lemma \ref{refereelemma}.
\end{proof}

\begin{corollary}
The nilpotent non-Lie Malcev algebra $\mathbb{M}$ is special: it is isomorphic
to a subalgebra of $A^-$ for some alternative algebra $A$.
\end{corollary}

\begin{proof}
The alternator ideal contains no elements of degree 1, and so the canonical map
from $\mathbb{M}$ to $A(\mathbb{M})$ is injective. Speciality of $\mathbb{M}$
also follows from Pchelintsev \cite{Pchelintsev}.
\end{proof}


\section*{acknowledgements}

We thank the referee for helpful comments; in particular, for pointing out the
need for (and providing a proof of) Lemma \ref{refereelemma}. This research was
partially supported by NSERC (Natural Sciences and Engineering Research Council
of Canada).



\begin{thebibliography}{99}

\bibitem
  {BremnerMurakamiShestakov}
  M. R. Bremner, L. I. Murakami and I. P. Shestakov:
  \emph{Nonassociative Algebras}.
  Chapter 69 of L. Hogben (editor),
  \emph{Handbook of Linear Algebra},
  Chapman \& Hall / CRC, Boca Raton, 2006.

\bibitem{Kuzmin1}
  E. N. Kuzmin:
  \emph{Malcev algebras and their representations}.
  Algebra Logic
  7 (1968) 48--69.

\bibitem
  {Kuzmin2}
  E. N. Kuzmin:
  \emph{Malcev algebras of dimension five over a field of characteristic zero}.
  Algebra Logic
  9 (1970) 416--421.

\bibitem{Malcev}
  A. I. Malcev:
  \emph{Analytic loops}.
  Mat{.} Sb{.} N.S.
  36 (1955) 569--576.

\bibitem{Pchelintsev}
  S. V. Pchelintsev:
  \emph{Speciality of metabelian Malcev algebras}.
  Math{.} Notes 74 (2003) 245--254.

\bibitem
  {PerezIzquierdoShestakov}
  J. M. P\'erez-Izquierdo and I. P. Shestakov:
  \emph{An envelope for Malcev algebras},
  J{.} Algebra
  272 (2004) 379--393.

\bibitem{Sagle}
  A. A. Sagle:
  \emph{Malcev algebras}.
  Trans{.} Amer{.} Math{.} Soc{.}
  101 (1961) 426--458.

\end{thebibliography}
\end{document}